\documentclass[12pt,reqno,a4paper]{article}
\pdfoutput=1 \widowpenalty10000 \clubpenalty10000

\usepackage{fullpage}
\usepackage[usenames]{color}
\usepackage{amsmath,amssymb,amsfonts,amsthm}
\usepackage{graphicx}
\usepackage{booktabs}
\usepackage{microtype}
\usepackage{graphics}
\usepackage{rotating}

\usepackage[colorlinks=true,linkcolor=webgreen,filecolor=webbrown,citecolor=webgreen]{hyperref}

\definecolor{webgreen}{rgb}{0,.5,0}
\definecolor{webbrown}{rgb}{.6,0,0}

\makeatletter\renewcommand{\fnum@figure}[1]{\small\figurename~\thefigure.}\makeatother
\makeatletter\renewcommand{\fnum@table}[1]{\small\tablename~\thetable.}\makeatother

\newcommand{\seqnum}[1]{\href{https://oeis.org/#1}{\rm \underline{#1}}}

\newcommand{\stirling}[2]{\genfrac{\{}{\}}{0pt}{}{#1}{#2}}

\DeclareMathOperator{\Li}{Li}

\makeatletter
\patchcmd{\@thm}{\trivlist}{\list{}{\leftmargin=2.5em \rightmargin=2.5em}}{}{}
\patchcmd{\@endtheorem}{\endtrivlist}{\endlist}{}{}
\makeatother

\theoremstyle{definition}
\newtheorem*{definition}{Definition}
\newtheorem*{conjecture}{Conjecture}
\newtheorem*{example}{Example}

\begin{document}

\begin{center}\vspace*{0.3cm}

{\LARGE\bf Cryptarithmically Unique Terms\smallskip\\ in Integer Sequences}

\vskip 1cm
\large
Dmytro~S.~Inosov\smallskip\smallskip\\
\normalsize
Institut f\"ur Festk\"orper- und Materialphysik\\
Technische Universit\"at Dresden\\
H\"ackelstra{\ss}e 3\\
01069 Dresden\\
Germany\smallskip\\
\href{mailto:dmytro.inosov@tu-dresden.de}{\tt dmytro.inosov@tu-dresden.de} \\

\vspace{.4cm}
\large
Emil Vlas\'ak\smallskip\smallskip\\
\normalsize
VVV-System, s.r.o.\\
V Podh\'aj\'i 776/30\\
40001 \'Ust\'i nad Labem\\
Czech Republic\smallskip\\
\href{mailto:emil@vlasak.biz}{\tt emil@vlasak.biz} \\
\end{center}

\vskip .2 in
\begin{abstract}
A cryptarithm (or alphametic) is a mathematical puzzle in which numbers are represented with words in such a way that identical letters stand for equal digits and distinct letters for unequal digits. An alphametic puzzle is usually given in the form of an equation that needs to be solved, such as \textsf{SEND}\,+\,\textsf{MORE}\,=\,\textsf{MONEY}. Alternatively, here we consider cryptarithms constrained not by an equation but by a particular subsequence of natural numbers, for example, perfect squares or primes. Such a cryptarithm has a unique solution if there is exactly one term in the sequence that has the corresponding pattern of digits. We call such terms cryptarithmically unique. Here we estimate the density of cryptarithmically unique terms using combinatorial considerations in an arbitrary sequence for which the overall density of terms is known. In particular, among all perfect squares below $10^{12}$, slightly less than one half are cryptarithmically unique, their density increasing toward larger numbers. Cryptarithmically unique prime numbers, however, are initially very scarce. Combinatorial estimates suggest that their density should drop below $10^{-300}$ for decimal lengths of approximately 1829 digits, but then it recovers and is asymptotic to unity for very large primes. Finally, we introduce and discuss primonumerophobic digit patterns that no prime number happens to have.
\end{abstract}

\section{\label{Sec:Introduction}Introduction}

\noindent Cryptarithms are well known in recreational mathematics as a type of a mathematical puzzle~\cite{VerbalArithmeticWiki}. The most famous and probably the oldest example, \textsf{SEND}\,+\,\textsf{MORE}\,=\,\textsf{MONEY}, was published in 1924 by Dudeney~\cite[pp.~97 and 214]{Dudeney24}. It has a unique solution, \mbox{\textsf{O}\,=\,0}, \mbox{\textsf{M}\,=\,1}, \mbox{\textsf{Y}\,=\,2}, \mbox{\textsf{E}\,=\,5}, \mbox{\textsf{N}\,=\,6}, \mbox{\textsf{D}\,=\,7}, \mbox{\textsf{R}\,=\,8}, and \mbox{\textsf{S}\,=\,9}. The rules of the game require that identical letters stand for equal digits, distinct letters for unequal digits, and none of the words start with a zero. The equation imposes a certain constraint on the possible values of letters that must guarantee a unique solution for the puzzle to make sense.

Alternatively, we can consider other types of constraints, for example, requiring that the number represented by the given word has a particular property or belongs to some predefined subset of natural numbers.
\begin{example}
$\sqrt{\text{\textsf{BIKINI}}}$ must be an integer.\smallskip\\
In this example, \textsf{BIKINI} must be a perfect square. For this particular word, the puzzle happens to have a unique solution, $\text{\textsf{BIKINI}}=190969=437^2$ or \textsf{B}\,=\,1, \textsf{N}\,=\,6, \textsf{I}\,=\,9, and \textsf{K}\,=\,0, whereas for some other words this may not be the case.
\end{example}\label{Example:Bikini}
If two integers allow for the same cryptarithmic representation, we call them \textit{cryptarithmically equivalent}. For example, 487242 is cryptarithmically equivalent to 531454 because they share the same pattern \textsf{``ABCDAD''} of repeating digits. It is obvious that cryptarithmically equivalent integers must be of the same length, meaning that to find all integers from a given subset that are equivalent to a given one, only a finite set of numbers needs to be checked. Every equivalence class is therefore finite, whereas the number of such equivalence classes can be infinite. Cryptarithmic equivalence classes are commonly referred to as \textit{digital types}~\cite{OEIS}.

Within a certain subset, those integers whose equivalence class is a singleton (contains no other elements except for this integer itself) are called \textit{cryptarithmically unique}. For example, the number 190969 from the above example is a cryptarithmically unique perfect square. It is natural to ask, how often cryptarithmically unique terms occur in certain integer sequences, for example, among squares or prime numbers. A related problem is to describe digital types that do not occur in the sequence. From purely combinatorial considerations, the higher the density of this sequence in $\mathbb{N}$, the less likely it is that for a particular integer none of its cryptarithmic equivalents belong to the sequence. It is also clear that integers with a small number of distinct digits are more likely to be cryptarithmically unique in a given subset, because the number of their cryptarithmic equivalents overall is much smaller.

\section{\label{Sec:CryptarithmicSignature}Cryptarithmic uniqueness}

\noindent To talk about cryptarithmic equivalence, it is useful to introduce a canonical form that would give the same value for equivalent integers and different values for nonequivalent ones.
\begin{definition}
For an arbitrary natural number $n\in\mathbb{N}$, its \textit{cryptarithmic signature}, $\text{CS}(n)$, is a number whose digits enumerate the decimal digits of the original number in order of their first occurrence. If the number has 10 distinct digits, the 10th digit is enumerated with zero.
\end{definition}
\noindent The function $\text{CS}(n)$ defined in this way is identical to the OEIS sequence \seqnum{A358497}($n$)~\cite{OEIS}. Here are several examples that illustrate this definition:
\begin{align*}
\text{CS}(237)&=123;\\
\text{CS}(190969)&=123242;\\
\text{CS}(2793297)&=1234132;\\
\text{CS}(8432916507)&=1234567890.
\end{align*}
Note that the cryptarithmic signature is an idempotent function, meaning that $\text{CS}(\text{CS}(n))=\text{CS}(n)$ for all $n\in\mathbb{N}$. By definition, the first digit of $\text{CS}(n)$ is always 1 and each successive digit to the right is at most one greater than the preceding digit (if we count 0 to represent 10), whereas the total number of digits $d$ is equal to that of the original number. 

In a similar way, one can define the cryptarithmic signature of an arbitrary word $W$ that has no more than 10 distinct letters. The solution of a cryptarithmic puzzle is therefore equivalent to solving the equation $\text{CS}(W)=\text{CS}(n)$ with respect to $n$ for some given word $W$ and with some constraint $n\in S\subset\mathbb{N}$.

We can now give more formal definitions of cryptarithmic equivalency and uniqueness using the cryptarithmic signature.
\begin{definition}
Two numbers $n\in\mathbb{N}$ and $m\in\mathbb{N}$ are called \textit{cryptarithmically equivalent} iff $\text{CS}(n)=\text{CS}(m)$.
\end{definition}
\begin{definition}
The number $n\in\mathbb{N}$ is called \textit{cryptarithmically unique} in some subset $S\subset\mathbb{N}$ if $n\in S$ and there exists no other number $m\in S$, $m\neq n$, such that $\text{CS}(n)=\text{CS}(m)$.
\end{definition}
\noindent These definitions are of course dependent on the base used. In particular, in base 2, every number is cryptarithmically unique because its first digit is always 1 and every other digit distinct from 1 must be 0, hence no two distinct numbers belong to the same cryptarithmic equivalence class. Most of the results in this paper will be given for the decimal base (base 10) but can be generalized to any arbitrary base $b$.

Note that the definition of cryptarithmic uniqueness is also applicable to any sequence $S=(s_i)_{i\,\geq\,0}$, $s_i\in\mathbb{N}$. The subset of cryptarithmically unique terms is independent of the order of terms in $S$, hence it can simply inherit the natural order of terms and be viewed as a subsequence of cryptarithmically unique terms. We let $\text{CU}(S)$ denote the subset (subsequence) of all cryptarithmically unique elements (terms) of $S$, given $S$ is a set (sequence).
\begin{definition}
$\text{CU}(S)=\{s_i\in S\,|\,\nexists\,j\neq i:\text{CS}(s_i)=\text{CS}(s_j),~s_j\in S\}$.
\end{definition}
\begin{definition}
A sequence for which $S=\text{CU}(S)$ is called a \emph{cryptarithmically unique sequence}.
\end{definition}
\begin{example}
The subsequence of natural numbers whose decimal expansion contains no digits larger than 1, which is \seqnum{A007088}$(n)$ (binary numbers)~\cite{OEIS} without the first term $n=0$, is cryptarithmically unique. Indeed, since the first digit of every term \seqnum{A007088}$(n)$, $n > 0$ is always 1, the terms are uniquely defined by the positions of zeros, hence no two terms can have the same digit pattern.
\end{example}

\noindent At most, a cryptarithmically unique sequence can have \seqnum{A164864}($d$)~\cite{OEIS} terms of decimal length $d$ (if it contains terms corresponding to every possible digit pattern).

\begin{definition}\noindent
A cryptarithmically unique sequence that has the maximal possible number of terms \seqnum{A164864}($d$) for every $d$ is called \emph{cryptarithmically complete}.
\end{definition}
\begin{example}
Trivial examples of such sequences are \seqnum{A358497}~\cite{OEIS} and \seqnum{A266946}~\cite{OEIS}, as they are cryptarithmically complete by definition.
\end{example}

\section{Cryptarithmic density}\label{Sec:CryptarithmicDensity}

If $S\subset\mathbb{N}$ is an infinite sequence, we can define its \textit{cryptarithmic density}:
\begin{definition}\noindent
For any integer sequence $S\subset\mathbb{N}$, the corresponding cryptarithmic density $D_S(d)$ is a function of the decimal length $d$ defined as the ratio of the number of cryptarithmically unique terms in $S$ of length $d$ to the total number of terms in $S$ of length $d$. For those $d$ for which no terms exist, $D_S(d)$ is defined to be 1:
\begin{displaymath}\label{Eq:DefinitionDS}
D_S(d)=\begin{cases}
       \,1, & \text{\hspace{-2mm}if $\{n\in S:\lfloor\log_{10}n\rfloor+1=d\}=\varnothing$;}\vspace{2pt}\\
       \displaystyle\,\frac{\#\{n\in\text{CU}(S):\lfloor\log_{10}n\rfloor+1=d\}}{\#\{n\in S:\lfloor\log_{10}n\rfloor+1=d\}}, & \text{\hspace{-2mm}otherwise.}
       \end{cases}
\end{displaymath}
\end{definition}
\noindent According to this definition, $D_S(d)\equiv1$ for any cryptarithmically unique sequence.
\begin{definition}\noindent
A sequence $S$ is called \emph{asymptotically cryptarithmically unique} if \mbox{$\lim_{d\rightarrow\infty} D_S(d) = 1$}.
\end{definition}
\begin{definition}\noindent
A sequence $S$ is called \emph{asymptotically cryptarithmically complete} if it is asymptotically cryptarithmically unique and, in addition,
\begin{equation}\label{Eq:DefinitionComplete}
\lim_{d\rightarrow\infty} \frac{\#\{n\in\text{CU}(S):\lfloor\log_{10}n\rfloor+1=d\}}{\text{\seqnum{A164864}}(d)}=1.
\end{equation}
\end{definition}

The asymptotic behavior of cryptarithmic density can be estimated from the asymptotic density of $S$ in $\mathbb{N}$ using combinatorial considerations under certain assumptions about the random distribution of decimal digits in $S$. Let us assume that we know the asymptotic probability $p(n)$ that a randomly chosen number $n$ belongs to $S$. For example, for perfect squares, this asymptotic probability is given by $1/\sqrt{n}$, while for prime numbers, it is asymptotic to $1/\ln(n)$ according to the prime number theorem~\cite[chapter~3, pp.~32--47]{Derbyshire04}. As we are interested in the average density of terms with a certain decimal length $d$, we can estimate it as
$$\tilde{p}(d)=10^d p(10^d) - 10^{d-1} p(10^{d-1}).$$

We start by considering a $d$-digit number \mbox{$n\in S$} that has exactly $k$ distinct digits $q_1, q_2,\ldots, q_k$ in its decimal representation, such that $q_i\neq q_j$ for $i\neq j$. The decimal representation of this number $n$ is given by $q_{i_1}q_{i_2}\ldots{\kern.5pt}q_{i_d}$, where $q_{i\!_j}\in\{q_1,\ldots,q_k\}$ for $j=1,\ldots,d$ and $q_{i_1}\neq0$. The total number of its cryptarithmic equivalents in $\mathbb{N}$ (including this number itself) is
$$N_\text{c}(k)=\frac{9\cdot9!}{(10-k)!}.$$
Note that $N_\text{c}(k)$ is independent of the length $d$ and only depends on the number of distinct digits $k$.

The cryptarithmic uniqueness of $n$ requires that none of its $N_\text{c}(k)-1$ cryptarithmic equivalents belong to $S$. The probability that this requirement is fulfilled can be estimated from combinatorial considerations under the assumption that these individual conditions are independent,
\begin{equation}\label{Eq:ProbOfUniqueness}
P(n,k)=\bigl(1-\tilde{p}(d)\bigr)^{N_\text{c}(k)-1}=\bigl(1-\tilde{p}(d)\bigr)^{\frac{9\cdot9!}{(10-k)!}-1},
\end{equation}
where $d=\lfloor\log_{10}n\rfloor+1$ is the integer length of $n$. Furthermore, we will need the probability that a randomly chosen $d$-digit number has exactly $k$ different digits, which is given by
\begin{equation}\label{Eq:ProbOfDistinctDigits}
p_k(d)=\frac{\stirling{d}{k}\binom{10}{k}k!}{10^d}=\frac{9!}{10^{d-1}(10-k)!}\sum_{i=1}^k\frac{(-1)^{k-i}i^{\,d}}{(k-i)!\,i!}.
\end{equation}
Here $\stirling{d}{k}$ is the Stirling number of the second kind, which gives the total number of distinct patterns in which $k$ different digits can be arranged into a number of length $d$, $\binom{10}{k}$ is the binomial coefficient that stands for the number of possible ways to assign decimal values to these $k$ different digits, and $k!$ gives the number of ways to order these $k$ different values. As expected, $\sum_{k=1}^{10} p_k(d)=1$ for any $d$.

\begin{figure}[t!]
\centerline{\includegraphics[width=0.8\linewidth]{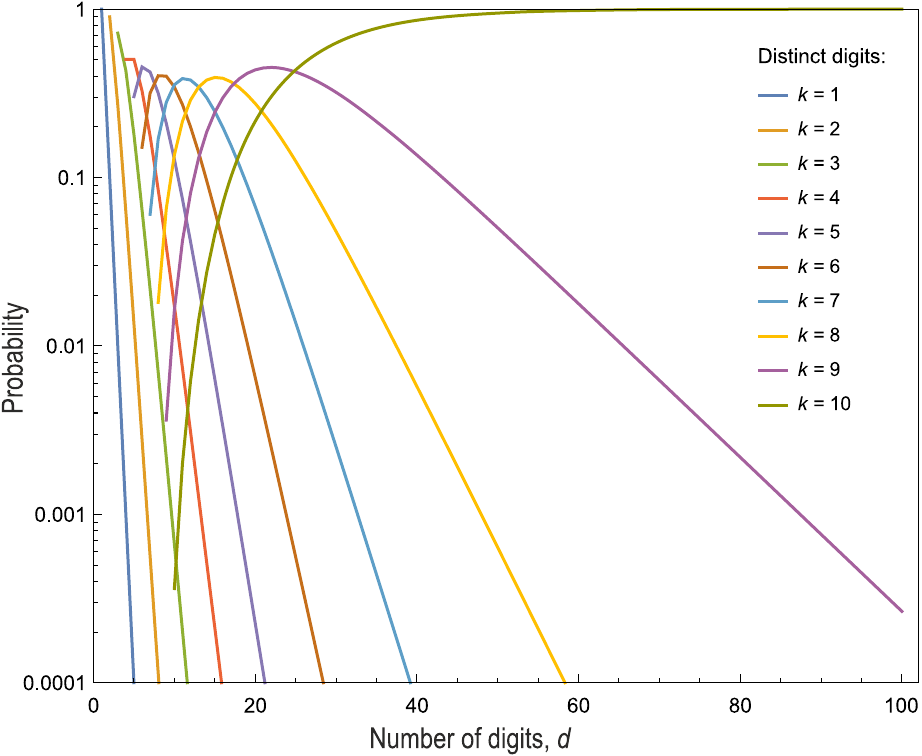}}
\caption{The probability $p_k(d)$ that a randomly chosen natural number of length $d$ has exactly $k$ distinct digits, calculated using Eq.~\eqref{Eq:ProbOfDistinctDigits}. All curves sum up to unity.}\label{Fig:Probability}
\end{figure}

The resulting probabilities are plotted in Figure~\ref{Fig:Probability}. One can see that for every $k$ there is a range of integer lengths in which integers with that number of distinct digits dominate. Finally, numbers in which all 10 distinct digits are present overwhelm at sufficiently large $d$.

To estimate the cryptarithmic density of an integer sequence, we calculate the combinatorial probability that a random element $n\in S$ of length $d=\lfloor\log_{10}n\rfloor+1$ is cryptarithmically unique. This is done by multiplying the probability that a certain term of length $d$ has exactly $k$ distinct digits, given by Eq.~\eqref{Eq:ProbOfDistinctDigits}, with the probability of cryptarithmic uniqueness, given by Eq.~\eqref{Eq:ProbOfUniqueness}, then summing up over all $k$:
\begin{equation}\label{Eq:CombinatorialUniquenessProb}
P(d)=\sum_{k=1}^{10}p_k(d)P(n,k)=\sum_{k=1}^{10}\frac{\stirling{d}{k}\binom{10}{k}k!}{10^d}\bigl(1-\tilde{p}(d)\bigr)^{\frac{9\cdot9!}{(10-k)!}-1}.
\end{equation}
Individual terms in this sum approximate the density of cryptarithmically unique terms of the sequence with exactly $k$ distinct digits, whereas the sum itself corresponds to the total cryptarithmic density.
To understand the asymptotic behavior of this function, it is useful to write out a few first terms explicitly:
\begin{align*}
P(d)&=\frac{1}{10^{d-1}}\bigl(1-\tilde{p}(d)\bigr)^{8} +\frac{9}{10^{d-1}}(2^{d-1}-1)\bigl(1-\tilde{p}(d)\bigr)^{80}\\
    &+\frac{36}{10^{d-1}}(3^{d-1}-2^d+1)\bigl(1-\tilde{p}(d)\bigr)^{647}+\,\cdots.
\end{align*}
We can see that after opening all the brackets and grouping linearly dependent terms, one obtains a sum of the following form:
\begin{equation}
P(d)=\sum_{l=0}^{9\cdot9!-1}\sum_{k=1}^{10}A_{kl}\Bigl(\frac{k}{10}\Bigr)^{\!d}\!\tilde{p}(d)^l,
\end{equation}
where $A_{kl}$ are constants. All terms in this sum are asymptotically negligible at $d\rightarrow\infty$ except for those with $k=10$. Furthermore, if we assume in addition that the distribution $\tilde{p}(d)$ satisfies $\lim_{d\rightarrow\infty} \tilde{p}(d)=0$, all terms with $l>0$ also become negligible with respect to the single term corresponding to $k=10$ and $l=0$, which is asymptotic to unity since $\lim_{d\rightarrow\infty}p_{10}(d)=1$. Therefore, for any such sequence, one expects that $\lim_{d\rightarrow\infty} P(d) = 1$. This should hold, in particular, for cryptarithmically unique primes with $\tilde{p}(d)\sim 1/(d\ln\!10)$ and squares with $\tilde{p}(d)\sim10^{-d/2}$.

While the following statement is certainly not true for every sequence, for certain ``well-behaved'' sequences it is natural to expect that $P(d)$ should asymptotically converge to the actual cryptarithmic density $D_S(d)$. On the one hand, the sequence $S$ must have an asymptotically vanishing density, $\lim_{d\rightarrow\infty} \tilde{p}(d)=0$, and on the other hand an increasing number of terms of every length $d$ so that $\lim_{d\rightarrow\infty} \#\{n\in S:\lfloor\log_{10}n\rfloor+1=d\}=\infty$, which guarantees that statistical fluctuations of the number of cryptarithmically unique terms vanish at large integer lengths. For such well-behaved sequences, asymptotic cryptarithmic uniqueness can be conjectured from combinatorial considerations by proving that $\lim_{d\rightarrow\infty} P(d) = 1$ and then hoping that it implies $\lim_{d\rightarrow\infty} D_S(d) = 1$.

\section{Cryptarithmically unique squares}

The number 190969, represented by the word \textsf{``BIKINI''} in the example on page \pageref{Example:Bikini}, is a cryptarithmically unique square, meaning that the equation $\text{CS}(n^2)=\text{CS}(190969)$ has a unique solution in natural numbers, $n=437$. This can be proven by simply checking all integers between 317 and 999 whose square has 6 digits. But how likely is it that a randomly chosen perfect square is cryptarithmically unique? And how does the density of cryptarithmically unique squares change as one goes to larger and larger numbers? Their sequence \seqnum{A374267}~\cite{OEIS} starts with $1444=38^2$, $7744=88^2$, $14884=122^2,\ldots$ The number of terms for different integer lengths is given in Table~\ref{Tab:NumOfUniques} along with the corresponding percentages with respect to the total number of perfect squares of the same length. These percentages describe the average density of cryptarithmically unique squares. For comparison, we also give the number of cryptarithmically unique triangular numbers, $n(n+1)/2$, in the same table. Because the densities of both sequences, $n^2$ and $n(n+1)/2$, are asymptotic to $1/\sqrt{n}$ for large $n$, we observe quite similar behavior in the densities of cryptarithmically unique terms of both sequences, which increase monotonically toward larger numbers (with the exception of two outliers for short triangular numbers where the number of hits is too low to draw any statistics from it).

\begin{table}[t!]
\begin{center}\small
\begin{tabular}{@{}lccccccccccc@{}}
    \toprule
        & \multicolumn{11}{c}{Number of decimal digits, $d$:}\\
        & 2 & 3 & 4 & 5 & 6 & 7 & 8 & 9 & 10 & 11 & 12 \\
	\midrule
squares: & 0 & 0 & 2 & 10 & 36 & 173 & 712 & 3522 & 17036 & 82953 & 386911\\
\%:      & 0 & 0 & 2.9 & 4.6 & 5.3 & 8.0 & 10.4 & 16.3 & 24.9 & 38.4 & 56.6\\
	\midrule
triangular: & 0 & 3 & 4 & 7 & 47 & 169 & 767 & 3682 & 18589 & 93389 & 460883\\
\%:      & 0 & 9.7 & 4.2 & 2.3 & 4.9 & 5.5 & 7.9 & 12.0 & 19.2 & 30.5 & 47.7\\
    \bottomrule
\end{tabular}
\caption{The total counts of cryptarithmically unique perfect squares~\cite{OEIS} and triangular numbers with the given number of decimal digits. The second row gives these numbers as percentages of the total number of squares and triangular numbers with the same number of digits, $d$.}\label{Tab:NumOfUniques}
\end{center}
\end{table}

\section{Cryptarithmically unique primes}

\subsection{Examples and general observations}

\begin{sidewaystable}\small
\begin{center}
\begin{tabular}{@{}l@{~}c@{~~}c@{~~}c@{~~}c@{~~}c@{~~}c@{~~}c@{~~}c@{}c@{}c@{}}
    \toprule
        & \multicolumn{10}{c@{}}{Number of decimal digits, $d$:\smallskip}\\
        & 2 & 7 & 8 & 9 & 10 & 11 & 12 & 13 & $\cdots$\hspace{-16pt} & 16\\
	\midrule
      & $k\!=\!1$: & $k\!=\!2$: & &         &            &            &               &               & & $k\!=\!2$:\\
      & 11 & 3333311 & 11818181 & 515115551 & 1110011101 & 11000010011 & 111515511511 & 1000001000111 & & 1000000010101111\\
      &    & 7771717 &          & 727722727 & 1161611161 & 11100000101 & 121212122221 & 1001001111001 & & 1000000101100111\\
      &    &         &          & 887887787 & 1411111441 & 11111313131 & 155511511151 & 1001100101011 & & 1000000110111001\\
      &    &         &          &           & 1411141411 & 11111313331 & 232322332223 & 1010000001001 & & 1000001011110001\\
      &    &         &          &           & 1717117117 & 11111661161 & 233233222333 & 1010000111101 & & 1000011110111011\\
      &    &         &          &           & 1911999919 & 11112121121 & 272277227777 & 1011010000001 & & 1000100000010101\\
      &    &         &          &           & 3311113111 & 11113131311 & 377733377737 & 1011110110111 & & 1000100110110001\\
      &    &         &          &           & 3313133311 & 11116611161 & 434334433343 & 1011111011101 & & 1000101010011101\\
      &    &         &          &           & 3333353533 & 11116666111 & 551551115111 & 1100000011111 & & 1000101011000011\\
      &    &         &          &           & 3333353533 & 11141414111 & 559595995559 & 1101001001101 & & 1000101111110111\\
      &    &         &          &           & 5151111551 & 11161116611 & 611616616111 & 1101011000101 & & \smash{\raisebox{-4pt}{$\vdots$}}\\
      &    &         &          &           & 5777777557 & 11177117111 & 661666661611 & 1110000001111 & & \\
      &    &         &          &           & 7373733337 & 11181118181 & 666111611611 & 1110110110111 & & $k\!=\!3$:\\
      &    &         &          &           & 7747447777 & 11181811111 & 667676667677 & 1111111010101 & & 1151135331533311\\
      &    &         &          &           & 7777111777 & 11616116111 & 676766677777 & 1111111110101 & & 1994191941199141\\
      &    &         &          &           & 8887788787 & 11711777771 & 677676767767 & 1111111111177 & & 2021201011100101\\
      &    &         &          &           & 9199119991 & 11717171177 & 737377377337 & 1111111122211 & & 2111252221521511\\
      &    &         &          &           & 9994449499 & 11717711177 & 766776677767 & 1111111141411 & & 3130101010001131\\
      &    &         &          &           &            & 11717777117 & 888898989889 & 1111111991999 & & 3133361113311661\\
      &    &         &          &           &            & 11991919991 & 911111999911 & 1111117177117 & & 3231113223111331\\
      &    &         &          &           &            & 11999199911 & 911199119911 & 1111119119191 & & 3331434141434341\\
      &    &         &          &           &            & 12121211111 & 991199199911 & 1111151515151 & & 3333333113136661\\
      &    &         &          &           &            & 13111311131 & 995555955599 & 1111161161111 & & 3773313111133711\\
      &    &         &          &           &            & \smash{\raisebox{-4pt}{$\vdots$}} & 999977997779 & \smash{\raisebox{-4pt}{$\vdots$}} & & \smash{\raisebox{-4pt}{$\vdots$}}\smallskip\smallskip\\
    \midrule
    Total: & 1 & 2 & 1 & 3 & 18 & 105 & 24 & 668 & & \hspace{-3pt}4235 ($k\!=\!2$) + 28 ($k\!=\!3$)\\
    \bottomrule
\end{tabular}
\caption{\small Cryptarithmically unique primes, \seqnum{A374238}~\cite{OEIS}, of different decimal lengths $d$ and with the number of distinct digits $k$. The columns with $d=3...\ldots6$ are skipped because no terms of these lengths exist. Note that lists with more than 24 terms are truncated. The bottom row gives the total number of cryptarithmically unique primes for these lengths, \seqnum{A376084}~\cite{OEIS}.}\label{Tab:ListOfUniquePrimes}
\end{center}
\end{sidewaystable}

It is easy to see that cryptarithmically unique primes with only a single distinct digit are prime repunits, i.e., prime numbers written as a sequence of 1's, \seqnum{A004022}~\cite{OEIS}. Indeed, any other number cryptarithmically equivalent to a repunit prime has the form $qqq{\kern.5pt}\cdots{\kern.5pt}qqq$ for some digit $q>1$, and is therefore divisible by $q>1$, hence nonprime. It is known that repunits with 2, 19, 23, 317, 1031, 49081, 86453, 109297, 270343, 5794777, and 8177207 digits are prime, \seqnum{A004023}~\cite{OEIS}. The smallest of them is 11, which is also the smallest cryptarithmically unique prime. Cryptarithmically unique primes with an arbitrary number of distinct digits $k$, \seqnum{A374238}~\cite{OEIS}, therefore represent a possible natural generalization (supersequence) of repunit primes, which then correspond to the special case of $k=1$.

Some cryptarithmically unique primes of different decimal lengths are listed in Table~\ref{Tab:ListOfUniquePrimes}. The number of such primes of different lengths $d$ for different numbers of distinct digits $k$ is given in Table~\ref{Tab:NumOfUniquePrimes}. There are no cryptarithmically unique primes of 1, 3, 4, 5, or 6 digits in length. We conjecture that except for the number 11, which is a repunit with only one distinct digit, all cryptarithmically unique primes of length smaller than 16 have only 2 distinct digits, the smallest of them is 3333311. The existence of terms $>10^{12}$ with 4 or more distinct digits has not been checked, but the combinatorial probability of their existence at such short integer lengths is vanishingly small, as we will show in the following. The total number of cryptarithmically unique primes of length $d$ is given by the OEIS sequence \seqnum{A376084}~\cite{OEIS}.

The smallest cryptarithmically unique prime with 3 distinct digits is 1151135331533311, it has $d=16$ decimal digits. There are in total 28 such primes of length 16 in addition to 4235 primes with 2 distinct digits. The first ten of each group are listed in Table~\ref{Tab:ListOfUniquePrimes}. Cryptarithmically unique primes with 4 or more distinct digits still await to be found.

\begin{table}[t!]
\begin{center}\small
\begin{tabular}{r@{~~}|@{~~}c@{~~~~}r@{~~~}r@{~~}|@{~}r}
    \toprule
     & \multicolumn{3}{c|@{~}}{$\hspace{-2ex}k$:} & \multicolumn{1}{c}{Total}\\
    \multicolumn{1}{c|@{~~}}{$d$} & 1 & 2 & 3 & \multicolumn{1}{c}{count}\\
	\midrule
    1 & 0 & 0 & 0 & 0 \\
    2 & 1 & 0 & 0 & 1 \\
    3 & 0 & 0 & 0 & 0 \\
    4 & 0 & 0 & 0 & 0 \\
    5 & 0 & 0 & 0 & 0 \\
    6 & 0 & 0 & 0 & 0 \\
    7 & 0 & 2 & 0 & 2 \\
    8 & 0 & 1 & 0 & 1 \\
    9 & 0 & 3 & 0 & 3 \\
    10 & 0 & 18 & 0 & 18 \\
    11 & 0 & 105 & 0 & 105 \\
    12 & 0 & 24 & 0 & 24 \\
    13 & 0 & 668 & 0 & 668 \\
    14 & 0 & 1129 & 0 & 1129 \\
    15 & 0 & 1306 & 0 & 1306 \\
    16 & 0 & 4235 & 28 & 4263 \\
    17 & 0 & 17100 & 220 & 17320 \\
    18 & 0 & 6704 & 30 & 6734 \\
    19 & 1 & 77057 & 4736 & 81794 \\
    20 & 0 & 110141 & 15834 & 125975 \\
    21 & 0 & 172639 & 7832 & 180471 \\
    22 & 0 & 499974 & 352605 & 852579 \\
    \bottomrule
\end{tabular}
\caption{Number of cryptarithmically unique primes of different lengths $d$ with different numbers of distinct digits $k$ and their total count, \seqnum{A376084}($d$)~\cite{OEIS}. For the total number of possible digit patterns that can result in primes, see Tables~\ref{Tab:DivisibilityRules} and \ref{Tab:PrimonumerophobicPatterns}.}\label{Tab:NumOfUniquePrimes}
\end{center}
\end{table}

Some of the cryptarithmically unique primes are simultaneously palindromes, \seqnum{A002113}~\cite{OEIS}. Because palindromicity is a property of the digit pattern itself, it applies to whole cryptarithmic equivalence classes rather than individual numbers, and therefore the sets of \emph{palindromic cryptarithmically unique primes} and \emph{cryptarithmically unique palindromic primes} coincide, hence these terms can be used synonymously. Cryptarithmically unique palindromic primes, \seqnum{A376118}~\cite{OEIS}, form a subsequence of both \seqnum{A374238}~\cite{OEIS} and \seqnum{A002385}~\cite{OEIS}. All palindromic primes except 11 have an odd number of digits, because any palindrome of even length is divisible by 11. The shortest cryptarithmically unique palindromic primes that are not repunits are the three 11-digit palindromes 11141414111, 11999199911, and 13111311131, followed by 12 palindromes of length 13, then 3 of length 15, then 73 of length 17, and so on. The shortest cryptarithmically unique palindromic primes with 3 distinct digits are the three 23-digit primes 11155511521212511555111, 11851585158585158515811, and 70009770799099707790007. Terms with 4 or more distinct digits still await to be found.

\subsection{Divisibility rules}

As one can see in Table~\ref{Tab:NumOfUniquePrimes}, the counts of cryptarithmically unique primes of different lengths $d$ do not form a monotonic sequence. Anomalously low values are observed, in particular, for $d=12$ and 18\,---\,lengths with a large number of divisors. This can be easily traced back to the divisibility rules.

\begin{sidewaystable}\small
\begin{tabular}{@{\,}r@{~~~}|@{~~}r@{~\,}r@{~\,}r@{~\,}r@{~\,}r@{~\,}r@{~\,}r@{~\,}r@{~\,}r@{~\,}r@{~~}|@{~\,}r@{\,}}
    \toprule
    $d\!\setminus\!k$\hspace{-0.6em} & 1 & 2 & 3 & 4 & 5 & 6 & 7 & 8 & 9 & 10 & Total\\
	\midrule
    1 & 1 & & & & & & & & & & 1\\
    2 & 1 & 1 & & & & & & & & & 2\\
    3 & 0 & 3 & 1 & & & & & & & & 4\\
    4 & 0 & 4 & 6 & 1 & & & & & & & 11\\
    5 & 0 & 15 & 25 & 10 & 1 & & & & & & 51\\
    6 & 0 & 12 & 84 & 65 & 15 & 1 & & & & & 177\\
    7 & 0 & 63 & 301 & 350 & 140 & 21 & 1 & & & & 876\\
    8 & 0 & 80 & 868 & 1672 & 1050 & 266 & 28 & 1 & & & 3965\\
    9 & 0 & 171 & 2745 & 7770 & 6951 & 2646 & 462 & 36 & 1 & & 20782\\
    10 & 0 & 370 & 8680 & 33505 & 42405 & 22827 & 5880 & 750 & 45 & 0 & 114462\\
    11 & 0 & 1023 & 28501 & 145750 & 246730 & 179487 & 63987 & 11880 & 1155 & 55 & 678568\\
    12 & 0 & 912 & 69792 & 583438 & 1373478 & 1322896 & 627396 & 159027 & 22275 & 1705 & 4160919\\
    13 & 0 & 3965 & 261495 & 2532517 & 7508501 & 9321312 & 5715424 & 1899612 & 359502 & 38610 & 27640938\\
    14 & 0 & 6412 & 745213 & 10230045 & 39913335 & 63383453 & 49324240 & 20912320 & 5135130 & 752752 & 190402900\\
    15 & 0 & 10200 & 2111715 & 40954050 & 209365520 & 420693273 & 408741333 & 216627840 & 67128490 & 12662650 & 1378295071\\
    16 & 0 & 22848 & 6703104 & 169164352 & 1091809304 & 2732339232 & 3281304992 & 2141723140 & 820784250 & 192453690 & 10436304912\\
    17 & 0 & 65535 & 21457825 & 694337290 & 5652751651 & 17505749898 & 25708104786 & 20415995028 & 9528822303 & 2758334150 & 82285618466\\
    18 & 0 & 55638 & 53630973 & 2658403776 & 28566917028 & 110396757606 & 197422820121 & 189029524401 & 106175032875 & 37112163803 & 671415306221\\
    19 & 1 & 262143 & 193448101 & 11259666950 & 147589284710 & 693081601779 & 1492924634839 & 1709751003480 & 1144614626805 & 472253194985 & 5671667723793\\
    20 & 0 & 399820 & 553780990 & \multicolumn{1}{@{}l}{\ldots} & & & & & & & \multicolumn{1}{@{}c@{\,}}{\ldots}\\
    \bottomrule
\end{tabular}
\caption{The number of digit patterns of length $d$ with $k$ distinct digits that satisfy no divisibility rules and the corresponding row sums, \seqnum{A376918}($d$)~\cite{OEIS}.}\label{Tab:DivisibilityRules}
\end{sidewaystable}

\begin{example}
Whenever the total number of digits is divisible by 3, certain digit patterns cannot result in primes. If every distinct digit occurs in the pattern a number of times that is divisible by 3, the sum of digits is also divisible by 3, and therefore the number cannot be prime. For example, for $n=12$ all patterns consisting of 2 distinct digits A and B that both occur a number of times divisible by 3 cannot produce primes and therefore do not contribute to the total count. This includes $\binom{12}{3}=220$ patterns with groups of 9 and 3 equal digits (such as \textsf{``AABABAAAABAA''} and alike) and $\frac{1}{2}\binom{12}{6}=462$ patterns with two groups of 6 equal digits (such as \textsf{``ABBABABAABBA''} and alike). The total number of 682 excluded patterns constitutes approximately one third of all $\stirling{12}{2}=2047$ possible 12-digit patterns with 2 distinct digits.
\end{example}

\begin{table}\small
\begin{center}
\begin{tabular}{@{\,}r@{~~~}|@{~~}r@{~\,}r@{~\,}r@{~\,}r@{~\,}r@{~\,}r@{~\,}r@{~\,}r@{~\,}r@{\,}}
    \toprule
    $d\!\setminus\!k$\hspace{-0.6em} & 1 & 2 & 3 & 4 & 5 & 6 & 7 & 8 & 9 \\
	\midrule
 2 & 1\\                                                                         
 3 & 1\\                                                                  
 4 & 1 &           3\\                                                       
 5 & 1 &           0\\                                             
 6 & 1 &          19 &         6\\    
 7 & 1 &           0 &         0\\  
 8 & 1 &          47 &        98 &        29\\ 
 9 & 1 &          84 &       280 &         0\\
10 & 1 &         141 &       650 &       600 &       120\\
11 & 1 &           0 &         0 &         0 &         0\\ 
12 & 1 &        1135 &     16734 &     28063 &      5922 &       756\\ 
13 & 1 &         130 &       130 &        13 &         0 &         0\\ 
14 & 1 &        1779 &     43757 &    161700 &    161700 &     52920 &     5040\\ 
15 & 1 &        6183 &    263386 &   1401900 &   1401400 &         0 &        0\\        
16 & 1 &        9919 &    438582 &   2634549 &   4381246 &   2587326 &   577612 &   40913\\                
17 & 1 &           0 &         0 &         0 &         0 &         0 &        0 &       0\\
18 & 1 &       75433 &  10808037 & 140403209 & 391178517 & 290493433 & 39663279 & 6540609 & 362880\\                 
19 & 1 &           0 &         0 &         0 &         0 &         0 &        0 &       0 &      0\\
20 & 1 &      124467 &  26825456 & 514583021 & \multicolumn{1}{@{}l}{\ldots}\\                              
21 & 1 &      369141 & 194219719 & \multicolumn{1}{@{}l}{\ldots}\\
22 & 1 &      429507 & 201337323 & \multicolumn{1}{@{}l}{\ldots}\\                                 
23 & 1 &           0 &         0 & 0 & 0 & 0 & 0 & 0 & 0 \\
24 & 1 &     4797007 & \multicolumn{1}{@{}l}{\ldots}\\                             
25 & 1 &      416965 & \multicolumn{1}{@{}l}{\ldots}\\                                   
26 & 1 &     6114993 & \multicolumn{1}{@{}l}{\ldots}\\                                   
27 & 1 &    22482399 & \multicolumn{1}{@{}l}{\ldots}\\                                   
28 & 1 &    28867895 & \multicolumn{1}{@{}l}{\ldots}\\                                                                            
29 & 1 &      111650 & \multicolumn{1}{@{}l}{\ldots}\\                                                                            
30 & 1 &   306153841 & \multicolumn{1}{@{}l}{\ldots}\\                                                                            
31 & 1 &      384524 & \multicolumn{1}{@{}l}{\ldots}\\                                                                            
32 & 1 &   507438239 & \multicolumn{1}{@{}l}{\ldots}\\                                                                            
33 & 1 &  1483501077 & \multicolumn{1}{@{}l}{\ldots}\\                                                                            
34 & 1 &  1242075013 & \multicolumn{1}{@{}l}{\ldots}\\                                                                            
35 & 1 &   743845628 & \multicolumn{1}{@{}l}{\ldots}\\                                                                            
36 & 1 & 19710473035 & \multicolumn{1}{@{}l}{\ldots}\\                                                                            
37 & 1 &       34299 & \multicolumn{1}{@{}l}{\ldots}\\                                                                            
38 & 1 & 17721793659 & \multicolumn{1}{@{}l}{\ldots}\\
39 & 1 & \multicolumn{1}{@{}c}{\ldots}\\
\bottomrule
\end{tabular}
\caption{The number of partitions $T(d,k)$ of the repunit \seqnum{A002275}($d$) into $k$ nonzero complementary binary vectors having a common divisor $> 1$ in base 10, see \seqnum{A378761}~\cite{OEIS}. All terms with $k > d/2$ are omitted as trivial zeros. Unknown terms are marked with ``$\ldots$''. In every line, a zero cannot be followed by nonzero terms. The sum of the first two columns ($k=1$ and $k=2$) is a separate sequence, \seqnum{A378511}($d$), and the row sums are given in \seqnum{A385539}~\cite{OEIS}.}\label{Tab:RepunitPartitions}
\end{center}
\end{table}

\begin{example}
Another subset of patterns that cannot produce primes is excluded because of the divisibility by 11. If the distinct digits \textsf{A}, \textsf{B} etc. occur on odd and even positions in the pattern the same number of times modulo 11, then the difference of the sums of digits on odd and even positions is divisible by 11 independently of the values of these digits, hence no number with this pattern can be prime. For example, the 26-digit pattern \textsf{``ABACACBCACABACBCACACACACAC''} has 11 \textsf{A}'s on odd positions, 2 \textsf{B}'s on odd and 2 \textsf{B}'s on even positions, and 11 \textsf{C}'s on even positions. Therefore the sums of digits on odd and even positions are $11\text{\textsf{A}}+2\kern.75pt\text{\textsf{B}}$ and $2\kern.75pt\text{\textsf{B}}+11\text{\textsf{C}}$, respectively, so their difference is divisible by 11 independently of the values of \textsf{A}, \textsf{B}, and \textsf{C}. Among 12-digit patterns with 2 distinct digits, 361 patterns can be excluded in this way, 100 of which have already been excluded due to the divisibility by 3. That adds 261 new patterns that cannot produce primes to the exclusion list.
\end{example}

There are two kinds of divisibility rules for digital types. The \emph{first kind}, illustrated by the two examples above, can be generalized to all prime divisors of $10^r-1$ or $10^r+1$, for which the following divisibility rules are known~\cite{ShaileshShirali}:
\begin{itemize}
\item If $p$ divides $10^r-1$, then the divisibility by $p$ can be tested by checking if the sum of $r$-digit blocks is divisible by $p$.
\item If $p$ divides $10^r+1$, then the divisibility by $p$ can be tested by checking if the alternating sum of $r$-digit blocks is divisible by $p$.
\end{itemize}
Using these rules, certain patterns for which $d$ is divisible by $r$ are excluded when considering the normal or alternating sum of $r$-digit blocks. In this way, 36 additional 12-digit patterns with 2 distinct digits are excluded by considering the alternating sum of two-digit blocks due to the divisibility by 101, then 144 more whenever the alternating sum of three-digit blocks is divisible by 7 or 13, and finally 12 more for which the alternating sum of six-digit blocks is divisible by 9901. In total, out of 2047 possible patterns of length 12 with $k=2$, only 912 remain that could in principle correspond to primes (see Table~\ref{Tab:DivisibilityRules}).

A general approach to find digital types that satisfy a divisibility rule of the first kind is to express the pattern with distinct digits \textsf{A}, \textsf{B},\,\textellipsis\ as a linear combination of binary vectors: $\text{\textsf{A}} \cdot X_1+\text{\textsf{B}}\cdot X_2 + \cdots$. The pattern coefficients $X_i$, $1\leq i \leq k$, consist of 0's and 1's, with 1's on positions of the corresponding digit in the pattern. They form a $k$-tuple of binary vectors that sum up to a repunit and are mutually complementary, meaning that for every position the digit 1 occurs on that position in exactly one of the vectors. Whenever the coefficients $X_1$, $X_2\ldots$ have a common divisor of a different digital type from the one we started from, the original pattern is composite for any values of the distinct digits and therefore cannot produce primes. Digital types with $k$ distinct digits that cannot generate primes because of the divisibility rules of the first kind therefore correspond to partitions of the repunit \seqnum{A002275}$(d)$ into $k$ nonzero complementary binary vectors having a common divisor $>1$. The number of such partitions $T(d,k)$ is given in Table~\ref{Tab:RepunitPartitions} (OEIS sequences \seqnum{A378511}, \seqnum{A378761}, and \seqnum{A385539}~\cite{OEIS}). This table reveals several interesting patterns.

\begin{enumerate}
\item $T(d,k) = 0$ for $k > d/2$ because in such partitions at least one term must contain just a single 1 (with all other digits zero) and is, therefore, a power of 10. Hence it cannot have nontrivial common divisors with the repunit \seqnum{A002275}$(d)$. Such terms are skipped in the table as trivial zeros.
\item $T(n,1)=1$ for $n \geq 2$ because the only partition that counts toward this term is the trivial one consisting of the repunit itself.
\item Whenever $d$ is a term in \seqnum{A385537}~\cite{OEIS}, the repunit \seqnum{A002275}$(d)$ is coprime with any other nonzero binary vector of the same length, and therefore all terms with $k>1$ vanish: $T(\text{\seqnum{A385537}}(n),k)=0$ for $k\geq2$ and any $n \geq 1$. In particular, this happens for indices of prime repunits, \seqnum{A004023}, which form a subsequence of \seqnum{A385537}~\cite{OEIS}. Indeed, existence of a common divisor $q>1$ for all binary vectors in the partition implies that their sum should also be divisible by $q$, which can only be fulfilled for the trivial partition consisting of the repunit itself ($k=1$) iff $d$ is a term in \seqnum{A385537}.
\begin{example}
The row $d=7$ contains no nonzero terms beyond $k=1$ despite the corresponding repunit \seqnum{A002275}$(7)$ being composite. This is because 1111111 is a semiprime repunit, \seqnum{A102782}~\cite{OEIS}, whose both prime divisors 239 and 4649 are large, and neither of them divides any binary number smaller than 1111111 itself. The same is observed for $d=11$ and 17 that are also terms in \seqnum{A046413} (indices of repunit semiprimes)~\cite{OEIS}. This example illustrates the general pattern that $T(d,k)$ shows anomalously low values for those $d$ for which the number of distinct prime factors of the repunit, \seqnum{A095370}~\cite{OEIS}, is small, and the smallest of these prime factors, \seqnum{A067063}$(d)$, is large.
\end{example}
\item In every row, a zero cannot be followed by nonzero terms. In other words, if $T(d,k)=0$ for some $d$ and $k$, then $T(d,m)=0$ also for every $m \geq k$. Indeed, if some $m$-tuple of binary vectors existed that is counted toward $T(d,m)$, then an $(m-1)$-tuple obtained by summing any two of its vectors while leaving others unchanged would be counted toward $T(d,m-1)$. By induction, this leads to $T(d,k)>0$, which is a contradiction.
\begin{example}
In the row $d=13$, the last nonzero term is $T(13,4)=13$, followed by two nontrivial zeros below the $k=d/2$ diagonal. The same is observed in the row $d=15$ beyond the 5th term.
\end{example}
\end{enumerate}

The \emph{second kind} of divisibility rules applies only to patterns with $k=10$ distinct digits whose length $d$ has the form $10m+3q$ for some integer $m \geq 1$ and $q \geq 0$. Since the sum of all digits from 0 to 9 equals 45, which is divisible by 9, any natural number in which all distinct digits from 0 to 9 have the same number of occurrences modulo 3 is also divisible by 9, hence the corresponding digital type is identically divisible by 9 and cannot produce primes.

Table~\ref{Tab:DivisibilityRules} gives the number of digit patterns of length $d$ with $k$ distinct digits that remain after excluding those that satisfy at least one divisibility rule, as well as their total count (OEIS sequence \seqnum{A376918}~\cite{OEIS}). This is the number of patterns that need to be checked when searching for cryptarithmically unique primes, since all patterns not counted toward \seqnum{A376918}($d$) cannot produce primes by definition. Note that the column corresponding to $k=2$ is nonmonotonic, with anomalously low values corresponding to integer length that have a large number of divisors, which explains corresponding anomalies in the number of cryptarithmically unique primes, \seqnum{A376084}($d$). 

\subsection{Cryptarithmic density of primes}

\begin{figure}[t!]
\centerline{\includegraphics[width=0.95\linewidth]{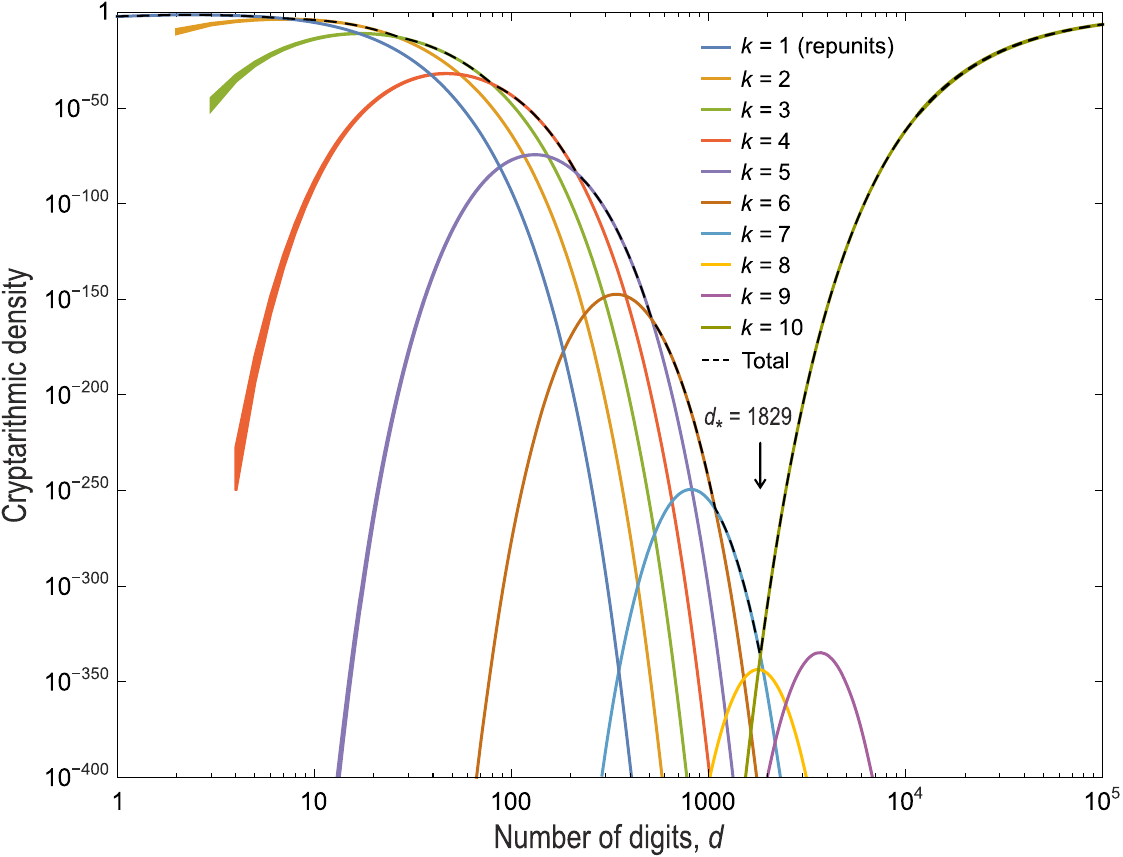}}
\caption{Combinatorial estimates of the density of cryptarithmically unique primes of length $d$ with exactly $k$ distinct digits, calculated with respect to the total number of primes estimated from the prime number theorem. The thickness of the lines represents the uncertainty range between the upper and lower estimates that result from two approximations of the prime-counting function, $x/\ln x$ and $\Li\,x$. The total density (sum of all curves) is plotted with the black dashed line.}\label{Fig:DensityPrimes}
\end{figure}

Now let us estimate the cryptarithmic density of primes using the results from section \ref{Sec:CryptarithmicDensity} and discuss its asymptotic behavior. Here we can employ two well-known asymptotic approximations for the prime-counting function $\pi(x)$:
\begin{align}
\pi(x)&\sim\frac{x}{\ln(x)};\\
\pi(x)&\sim\Li(x)=\int_2^x\frac{1}{\ln(t)}{\rm d}t.
\end{align}
These give us two expressions for the lower and upper bounds of the density of primes of length $d$, which is calculated as $(\pi(10^d)-\pi(10^{d-1}))/(9\cdot10^{d-1})$:
\begin{equation}
\tilde{p}_\text{l}(d)\sim\frac{1}{d\ln 10};~~\tilde{p}_\text{u}(d)\sim\frac{\Li(10^d)-\Li(10^{d-1})}{9\cdot10^{d-1}}.
\end{equation}
After substituting these into Eq.~\eqref{Eq:CombinatorialUniquenessProb}, we obtain the combinatorial estimates of the density of cryptarithmically unique primes of length $d$, which are plotted in Figure~\ref{Fig:DensityPrimes} separately for every number of distinct digits $k$ and as a sum (black dashed line). The interval between the lower and upper bounds corresponding to $\tilde{p}_\text{l}(d)$ and $\tilde{p}_\text{u}(d)$ is represented by the line thickness. Note that the divisibility rules may result in additional factors that are not taken into account in our calculation, as they would be barely visible on our logarithmic scale. The influence of the divisibility rules is expected to vanish asymptotically as $d\rightarrow\infty$.

The remarkable aspect of the resulting function is that it has a deep minimum at $d_\ast=1829$ where it equals approximately $4.0\cdot\!10^{-337}$. At $d_\ast$, cryptarithmically unique primes with 7 and 10 distinct digits are equally likely. However, as soon as cryptarithmically unique primes with $k=10$ overwhelm, the function recovers and asymptotically approaches unity, as expected.
\begin{conjecture}\noindent
The sequence of primes is asymptotically cryptarithmically unique.
\end{conjecture}
\noindent Indeed, the asymptotics $\tilde{p}(d)\sim 1/(d\ln\!10)$ in Eq.~\eqref{Eq:CombinatorialUniquenessProb} suggests that \mbox{$\lim_{d\rightarrow\infty} D_{\mathbb{P}}(d) = 1$}, under the assumption that large primes do not cluster into groups with the same digit patterns. This is a statement that would be impossible to reach without the asymptotic analysis of combinatorial probabilities, since actually calculating and counting cryptarithmically unique primes at such large lengths is computationally intractable. Any empirical observation based on the monotonically decreasing cryptarithmic density of primes with less than 1828 digits, for which cryptarithmic uniqueness can be realistically checked, would be misleading.
\begin{example}\noindent
The largest prime number known to date, which is the Mersenne prime \mbox{$M_\text{136,279,841}=2^{136,279,841}-1$}, has 41,024,320 decimal digits~\cite{Mersenne}, which is well beyond the right-hand cutoff of our plot in Figure~\ref{Fig:DensityPrimes}. At this length, the combinatorial probability that $M_\text{136,279,841}$ is cryptarithmically unique is approximately~96.6\%.
\end{example}

\begin{figure}[t!]
\centerline{\includegraphics[width=0.9\linewidth]{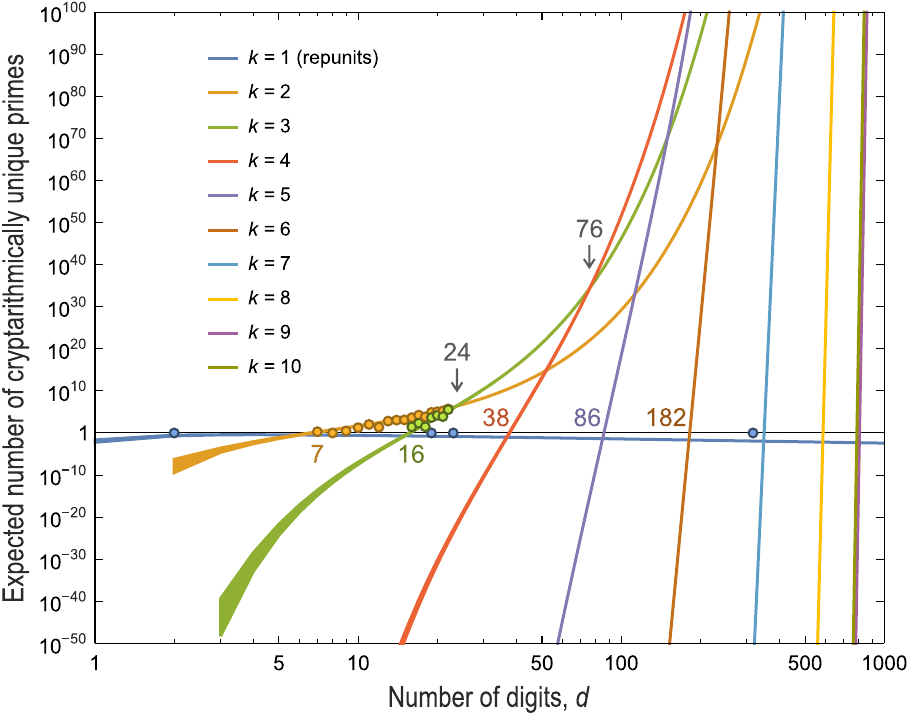}}
\caption{Combinatorial estimates of the expected number of cryptarithmically unique primes of length $d$ with exactly $k$ distinct digits (solid lines). The thickness of the lines represents the uncertainty range between the upper and lower estimates that result from two approximations of the prime-counting function, $x/\ln x$ and $\Li\,x$. Circles show the known actual counts of cryptarithmically unique primes from Table~\ref{Tab:NumOfUniquePrimes}.}\label{Fig:NumberOfPrimes}
\end{figure}

In Figure~\ref{Fig:NumberOfPrimes}, we also plot the expected number of cryptarithmically unique primes of different lengths with different number of unique digits. Circles show the actually observed number of primes from Table~\ref{Tab:NumOfUniquePrimes}. This plot shows, in particular, at which integer lengths the smallest primes with a given $k$ are to be expected. The length $d$ at which the expected number first exceeds 1 are written next to the curves. For $k=2$ and 3, the estimates are 7 and 16 digits, which exactly matches the actual lengths of the smallest cryptarithmically unique primes with 2 and 3 distinct digits in Table~\ref{Tab:NumOfUniquePrimes}. The expected length of the smallest cryptarithmically unique prime with 4 digits is 38, so finding it would be much more challenging.

We also see that every $k$ has a range of lengths in which primes with $k$ distinct digits dominate, for instance, cryptarithmically unique primes with 3 distinct digits are expected to become more common than those with 2 distinct digits at $d\geq24$, then primes with 4 distinct digits would take over at $d\geq76$, and so on.

\section{Primonumerophobic digit patterns}

As we have seen before, certain digit patterns (digital types) can never correspond to primes because of divisibility rules. However, some of the remaining patterns may still contain no primes ``accidentally'', without any connection to divisibility rules. This requires that all members of the cryptarithmic equivalence class represented by this pattern are nonprime. We call such patterns \emph{primo\-nu\-me\-ro\-phobic}.
\begin{definition}
A digit pattern (digital type) $W$ is called \emph{primonumerophobic} if it satisfies no divisibility rules (counts toward \seqnum{A376918}~\cite{OEIS}) and
$$\{n\in \mathbb{N}:\text{CS}(n)=\text{CS}(W)\}\cap\mathbb{P}=\varnothing.$$
\end{definition}
The existence of primonumerophobic digit patterns explains why the sequence \seqnum{A376918}($d$)~\cite{OEIS} does not coincide with \seqnum{A267013}($d$), which gives the actual number of distinct digit patterns corresponding to primes of length~$d$~\cite{OEIS}.

It is interesting to see how many primonumerophobic patterns exist and how their number depends on the integer length. Note that according to our definition, those patterns that can be proven to be nonprime based on divisibility rules (without the need to check every individual member) are not counted as primonumerophobic. In particular, no repdigit pattern ($k=1$) can be primonumerophobic, since any such number has the form $A(10^r-1)/9$, but for all prime factors of $(10^r-1)/9$ there is a divisibility rule based on the sum of $r$-digit blocks.

The number of primonumerophobic patterns of different lengths $d$ for $k=2$ and 3 is given in Table~\ref{Tab:NumOfUniquePrimes} together with the total number of available patterns that remain after divisibility-rule exclusions. The shortest patterns of this kind have length 10, these are \textsf{``AAABBBAAAB''}, \textsf{``AABABBBBBA''}, and \textsf{``ABAAAAABBB''}. There are in total 81 numbers corresponding to each pattern, none of which are prime. The shortest primonumerophobic patterns with 3 distinct digits have length 16, these are \textsf{``AAABBCABCCCAACCB''}, \textsf{``AABCBACAACCABCCC''}, \textsf{``ABBBACCBBCAAAABC''}, and \textsf{``ABBBCCBCABABBBBC''}. Here 648 different numbers corresponding to each pattern simultaneously happen to be nonprime.

\begin{table}[t!]
\begin{center}\small
\begin{tabular}{r@{~~}|@{~~}r@{~\,out of\,~}r@{~~}|@{~~}r@{~\,out of\,~}r@{~~}|@{~}r}
    \toprule
     & \multicolumn{4}{c@{~~}|@{~}}{$k$:} & \multicolumn{1}{c}{Total}\\
    \multicolumn{1}{c|@{~~}}{$d$} & \multicolumn{2}{c@{~~}|@{~~}}{2} & \multicolumn{2}{c@{~~}|@{~}}{3} & \multicolumn{1}{c}{count}\\
	\midrule
    1 & 0 & 0 & 0 & 0 & 0\\
    2 & 0 & 1 & 0 & 0 & 0\\
    3 & 0 & 3 & 0 & 1 & 0\\
    4 & 0 & 4 & 0 & 6 & 0\\
    5 & 0 & 15 & 0 & 25 & 0\\
    6 & 0 & 12 & 0 & 84 & 0\\
    7 & 0 & 63 & 0 & 301 & 0\\
    8 & 0 & 80 & 0 & 868 & 0\\
    9 & 0 & 171 & 0 & 2745 & 0\\
    10 & 3 & 370 & 0 & 8680 & 3\\
    11 & 32 & 1023 & 0 & 28501 & 32\\
    12 & 9 & 912 & 0 & 69792 & 9\\
    13 & 207 & 3965 & 0 & 261495 & 207\\
    14 & 362 & 6411 & 0 & 745209 & 362\\
    15 & 363 & 10046 & 0 & 2107764 & 363\\
    16 & 1444 & 22848 & 4 & 6703102 & 1448\\
    17 & 7583 & 65535 & 28 & 21457825 & 7611\\
    18 & 1837 & 55638 & \ldots & 53630973 & \multicolumn{1}{r}{\ldots}\\
    19 & 39137 & 262143 & \ldots & 193448101 & \\
    20 & 52480 & 399820 & \ldots & 553780990 & \\
    \bottomrule
\end{tabular}
\caption{Number of primonumerophobic digit patterns of different lengths $d$ with different numbers of distinct digits $k$. The second value (after ``out of'') is the total number of candidate patterns for the same $d$ and $k$ that remain after those patterns that satisfy any of the divisibility rules (i.e., can be proven nonprime without checking every individual number) have been excluded, according to Table~\ref{Tab:DivisibilityRules}. The last column gives the total number of primonumerophobic digit patterns, \seqnum{A377727}~\cite{OEIS}.}\label{Tab:PrimonumerophobicPatterns}
\end{center}
\end{table}
 
Palindromic patterns with an even number of digits cannot be primonumerophobic because they are divisible by 11, which can be checked using a divisibility rule based on the alternating sum of digits. However, \emph{primonumerophobic palindromes} with an odd number of digits do exist, and the shortest such palindromes are nine 13-digit patterns \textsf{``AAAABABABAAAA''}, \textsf{``AAAABBABBAAAA''}, \textsf{``AABAAABAAABAA''}, \textsf{``AABABBABBABAA''}, \textsf{``ABABABABABABA''}, \textsf{``ABABBBBBBBABA''}, \textsf{``ABBAABBBAABBA''}, \textsf{``ABBBAAAAABBBA''}, and \textsf{``ABBBABBBABBBA''}. The shortest primonumerophobic palindromes with 3 distinct digits are the three 25-digit patterns \textsf{``AABCCCABBCAACAACBBACCCBAA''}, \textsf{``ABABABACACCCCCCCACABABABA''}, and \textsf{``ABBCBCAAACAAAAACAAACBCBBA''}.

One can see that the number of primonumerophobic digit patterns also exhibits a local minimum at $d=12$, which has many divisors so that more than half of all possible patterns are excluded by divisibility rules. Subtracting the left number from the right one in a column corresponding to a particular $k$ in Table~\ref{Tab:PrimonumerophobicPatterns} gives the actual number of different existing prime patterns that are indeed realized for those $d$ and $k$ and contribute to \seqnum{A267013}($d$)~\cite{OEIS}. This allows calculating new conjectured terms in the sequence \seqnum{A267013} with less effort than with direct brute-force counting, assuming that contributions with $k>3$ decimal digits are zero as the combinatorial probability of their existence remains vanishingly small until much larger integer lengths.

From combinatorial considerations, one expects the number of primonumerophobic digit patterns to remain smaller than that of cryptarithmically unique primes, as one extra condition needs to be fulfilled: Now all $N_\text{c}(k)$ cryptarithmic equivalents must be nonprime rather than all but one. Apart from this one additional factor of $1-1/(d\ln10)$, which is asymptotic to 1, we expect these numbers to follow similar trends to those shown in Figs.~\ref{Fig:DensityPrimes} and \ref{Fig:NumberOfPrimes}. Therefore, it can be expected that primonumerophobic patterns as well as primonumerophobic palindromes with any given number of distinct digits $1 < k \leq 10$ do exist, even if those with 4 or more distinct digits still await to be found.

\subsection{Density of prime-generating digit patterns}

While we anticipate that the sequence of primes must be asymptotically cryptarithmically unique, it is easy to see that it is not asymptotically cryptarithmically complete. Indeed, the expected number of cryptarithmically unique primes of length $d$ is asymptotically the same as the total number of primes of the same length, which is asymptotic to $9\cdot10^{d-1}/(d\cdot\ln10)$. On the other hand, the total number of available digit patterns \seqnum{A164864}($d$) in the denominator of Eq.~\eqref{Eq:DefinitionComplete} is given by
$$
\text{\seqnum{A164864}}(d)=\sum_{k=1}^{10}\frac{1}{k!}\sum_{i=0}^k(-1)^{k-i}\binom{k}{i}i^d,
$$
which as a function of $d$ represents a sum of exponents with the leading term corresponding to $k=10$ and $i=10$. It is therefore asymptotic to $10^d/10!$ and results in
$$
\lim_{d\rightarrow\infty} \frac{\#\{n\in\text{CU}(\mathbb{P}):\lfloor\log_{10}n\rfloor+1=d\}}{\text{\seqnum{A164864}}(d)}=\frac{9\cdot9!}{\ln10}\lim_{d\rightarrow\infty} \frac{1}{d}=0.
$$
This is consistent with our expectation that almost all digit patterns become primonumerophobic as $d\rightarrow\infty$.

The fraction $F(d)$ of all digit patterns that long prime numbers are expected to exhibit on average (neglecting the divisibility-rule exclusions) among all the \seqnum{A164864}($d$) available patterns can be estimated using similar combinatorial considerations to those employed in the derivation of Eq.~\eqref{Eq:CombinatorialUniquenessProb}. The probability that a given pattern is not primonumerophobic requires that at least one of the $9\cdot9!/(10-k)!$ numbers corresponding to this pattern is prime, which we approximate with $1-\bigl(1-1/(d\ln10)\bigr)^{9\cdot9!/(10-k)!}$. After multiplying it with the total number of available patterns with $k$ distinct digits, $\stirling{d}{k}$, and summing over all $k$, we obtain the expected number of non-primonumerophobic patterns of length $d$. Then after dividing it by the total number of patterns \seqnum{A164864}($d$), we obtain the final result:
\begin{equation}\label{Eq:CombinatorialPatternsPrime}
F(d)\sim1-\sum_{k=1}^{10}\frac{\displaystyle\left(\!1-\frac{1}{d\ln10}\right)^{\frac{9\cdot9!}{(10-k)!}}\sum_{i=1}^{k}\frac{(-1)^{k-i}i^{\,d}}{(k-i)!i!}}
{\displaystyle\sum_{k=1}^{10}\sum_{i=1}^{k}\frac{(-1)^{k-i}i^{\,d}}{(k-i)!i!}}.
\end{equation}
In the limit $d\rightarrow\infty$, all terms in the numerator and denominator with $k<10$ and $i<10$ are asymptotically negligible, while the remaining leading term has the asymptotics \mbox{$F(d)\sim9\cdot9!/(d\ln10)$}, which we obtained above. This lets us conjecture that the ratio of the sequences \seqnum{A267013}~\cite{OEIS} and \seqnum{A164864}~\cite{OEIS} should have the same asymptotics:
\begin{conjecture}\noindent
The ratio $\displaystyle \frac{\text{\seqnum{A267013}}(d)}{\text{\seqnum{A164864}}(d)}$ is asymptotic to $\displaystyle\frac{9\cdot9!}{d\ln10}$ as $d\rightarrow\infty$.
\end{conjecture}
Note that the sum in Eq.~\eqref{Eq:CombinatorialPatternsPrime} is identical to that in Eq.~\eqref{Eq:CombinatorialUniquenessProb} after substituting $\tilde{p}(d)=\tilde{p}_1(d)$, except for one additional factor of $\bigl(1-1/(d\ln10)\bigr)$:
$$
F(d)=1-\left(\!1-\frac{1}{d\ln10}\right)P(d).
$$
The function $F(d)$ has a broad maximum at $d_\ast=1829$, where it reaches a value of approximately $1-4.0\cdot10^{-337}$, and then starts monotonically decreasing to zero, ultimately approaching the $9\cdot9!/(d\ln10)$ asymptote for much larger $d$. This means that 1829-digit patterns have the highest probability of containing primes in their equivalence classes, whereas the probability of finding cryptarithmically unique primes or primonumerophobic patterns at this length is minimal.

\bigskip

\hrule\bigskip
\noindent
2020 \emph{Mathematics Subject Classification:}
Primary: 00A08; Secondary 11A41, 05A05, 05A16.\\
\emph{Keywords:} digital type, digit pattern, cryptarithm, cryptarithmic uniqueness, prime, asymptotic enumeration, primonumerophobic digit pattern.
\bigskip
\hrule\bigskip
\noindent
(Concerned with sequences \seqnum{A267013}, \seqnum{A374238}, \seqnum{A374267}, \seqnum{A376084}, \seqnum{A376118}, \seqnum{A376918}, \seqnum{A377727}, \seqnum{A378154}, \seqnum{A378199}, \seqnum{A378511}, \seqnum{A378514}, \seqnum{A378761}, \seqnum{A385537}, \seqnum{A385539})
\bigskip
\hrule

\end{document}